\documentclass[11pt]{amsart}
\usepackage{amssymb,latexsym,amsmath,amsthm,amsfonts,amsxtra}

\theoremstyle{plain}
\newtheorem{theorem}{Theorem}
\newtheorem{corollary}{Corollary}

\newtheorem{lemma}{Lemma}

\newtheorem{remark}{Remark}

\theoremstyle{definition}
\newtheorem{definition}{Definition}

\theoremstyle{remark}

\numberwithin{equation}{section}

\numberwithin{equation}{section}

\newcommand{\ZZ}{\mathbb{Z}}
\newcommand{\Z}{\mathcal{Z}}
\newcommand{\K}{\mathbb{K}}

\newcommand{\RR}{\mathbb{R}}
\newcommand{\CC}{\mathbb{C}}

\newcommand{\B}{{\mathcal{B}}}
\newcommand{\BB}{{\mathbb{B}}}
\newcommand{\D}{{\mathbb{D}}}
\newcommand{\bfd}{{\mathbf{d}}}
\newcommand{\bfn}{{\mathbf{n}}}

\author{Manuel Lladser}

\title[Mixed Powers of Generating Functions]{Mixed powers of generating functions}


\begin{document}

\maketitle

\begin{abstract}
Given an integer $m\ge1$, let $\|\!\cdot\!\|$ be a norm in $\mathbb{R}^{m+1}$ and let $\mathbb{S}_+^{m}$ denote the set of points $\mathbf{d}=(d_0,\ldots,d_m)$ in $\mathbb{R}^{m+1}$ with nonnegative coordinates and such that $\|\mathbf{d}\|=1$.
Consider for each $1\le j\le m$ a function $f_j(z)$ that is analytic in an open neighborhood of the point $z=0$ in the complex plane and with possibly negative Taylor coefficients. Given $\mathbf{n}=(n_0,\ldots,n_m)$ in $\ZZ^{m+1}$ with nonnegative coordinates, we develop a method to systematically associate a parameter-varying integral to study the asymptotic behavior of the coefficient of $z^{n_0}$ of the Taylor series of $\prod_{j=1}^m \{f_j(z)\}^{n_j}$, as $\|\bfn\|\to\infty$. The associated parameter-varying integral has a phase term with well specified properties that make the asymptotic analysis of the integral amenable to saddle-point methods: for many $\mathbf{d}\in\mathbb{S}_+^m$, these methods ensure uniform asymptotic expansions for $[z^{n_0}]\prod_{j=1}^m \{f_j(z)\}^{n_j}$ provided that $\mathbf{n}/\|\mathbf{n}\|$ stays sufficiently close to $\mathbf{d}$ as  $\|\mathbf{n}\|\to\infty$. Our method finds applications in studying the asymptotic behavior of the coefficients of a certain multivariable generating functions as well as in problems related to the Lagrange inversion formula for instance in the context random planar maps.
\end{abstract}

\section{Introduction and Main Definitions}
\label{sec:introduction}

We start by introducing some notation that will be used consistently throughout this manuscript. In what follows, $m\ge1$ is a fixed integer. Define $\ZZ_+:=\{0,1,2,\ldots\}$ and $\RR_+:=[0,\infty)$. We use bold-face notation to refer to vectors in $\RR^{m+1}$ and denote the coordinates of a vector $\bfd$ as $(d_0,\ldots,d_m)$. The symbol $\bfn$ is reserved for elements in $\ZZ_+^{m+1}$. We let $\|\!\cdot\!\|$ be an arbitrary yet fixed norm in $\RR^{m+1}$. Define $\mathbb{S}_+^m:=\{\bfd\in\RR_+^{m+1}:\|\bfd\|=1\}$ and, for $\epsilon>0$, $\BB_\epsilon:=\{\bfd:\|\bfd\|<\epsilon\}$. The boundary of a set $\D\subset\mathbb{S}_+^m$ is denoted as $\partial\D$.

Given $\epsilon>0$ we write $[|z|<\epsilon]$ to refer to the open disk of radius $\epsilon$ centered at the origin in $\CC$.

For the remainder of this manuscript, $r>0$ is a fixed radius and, for each $1\le j\le m$, $f_j(z)$ is an analytic function over the disk $[|z|<r]$. We write $f_j$ as a shortcut for $f_j(z)$. For an arbitrary function $f(z)$ that is analytic about the origin we write $[z^n]\,f(z)$ to refer to the coefficient of $z^n$ of the Taylor series of $f(z)$ about $z=0$.

We are interested in the behavior of the Taylor coefficients of $\prod_{j=1}^mf_j^{n_j}$ for non-negative integer exponents $n_j$. More specifically, given $\bfn=(n_0,\ldots,n_m)$, we would like to systematically provide an asymptotic expansion for the coefficient
\begin{equation}
\label{ide:problem definition}
[z^{n_0}]\,\prod_{j=1}^mf_j^{n_j}
\end{equation}
as $\|\bfn\|\to\infty$. Since in a finite dimensional space all norms are equivalent, the condition $\|\bfn\|\to\infty$ is equivalent to the condition $\|\bfn\|_\infty:=\max\{|n_0|,\ldots,|n_m|\}\to\infty$. The terminology of {\it mixed powers of generating functions} used in the title of this manuscript aims to emphasize the following fact:  in order for the asymptotic analysis of the coefficients in (\ref{ide:problem definition}) to fall in the context of our discussion it suffices that at least one of the exponents $n_0,\ldots,n_m$ blows up to infinity. Although the actual norm used in $\RR^{m+1}$ is not theoretically relevant, it is worth to stress out that a suitable choice for $\|\!\cdot\!\|$ may considerably simplify the numerical analysis in an specific example. We take a considerably advantage of this fact in the applications discussed in sections~\ref{sec:application} and \ref{sec:planar maps} ahead.

Various special instances of (\ref{ide:problem definition}) can be found in the literature. Coefficients of this form include the well-studied case of the coefficients $[z^j]\,f(z)$ as $j\to\infty$~\cite{Wil,FlaSed}. It also includes the case $[z^j]\,f(z)^k$ as $j\to\infty$ or $k\to\infty$. Coefficients of this form occur frequently in Discrete Probability. For instance, if $X$ is a non-negative integer-valued random variable with moment generating function $f(z)=E(z^X)$, $[z^j]\,f(z)^k$ corresponds to the probability that $X_1+\ldots+X_k=j$, where $X_1,\ldots,X_k$ are independent copies of $X$. On the other hand, at the very core of Computer Science, coefficients of this form also occur in standard multinomial allocation problems (also called occupancy problems). For example, using the identity
\[\sum_{j_1,\ldots,j_k}{j\choose j_1,\ldots,j_k}=j!\cdot[z^j](e^z-1)^k\]
where the summation indices are such that $j_i\ge1$ and $j_1+\ldots+j_k=j$, it is not difficult to see that $j!{l\choose k}[z^j](e^z-1)^j$ represents the number of ways of allocating $j$ distinguishable balls to exactly $k$ of $l$ available urns of unbounded capacities~\cite{JohKot}.

Terms of the form of (\ref{ide:problem definition}) also include the coefficients $[z^j]\,f(z)^kg(z)^l$. Coefficients of this last form, with $k/j\to1$ and $l$ bounded, arise in the context of the Lagrange inversion formula~\cite{GouJac}: if $f(0)\ne0$ and $h(t)$ is defined implicitly through the relation $h(t)=t\cdot f(h(t))$ then
\[[t^j]\,g(h(t))=\frac{1}{j}\cdot[z^{j-1}]\,\{f(z)\}^jg'(z).\]

All the above instances of~(\ref{ide:problem definition}) have the following in common: Cauchy's formula is used to represent the coefficients in~(\ref{ide:problem definition}) via a contour-integral which is then analyzed using {\it saddle-point methods} (see~\cite{BleHan} for a comprehensive discussion of this method as well as other related methods). This method has two variants: the non-coalescing, and the coalescing version. While the non-coalescing version has been used systematically in the context of asymptotic enumeration, it was not until recently that Banderier et al~\cite{BFSS_b} used the more elaborated {\it coalescing saddle-point method} to unravel the asymptotic distribution of the core size of random planar maps. Their analysis relies on studying the asymptotic behavior of certain coefficients of the form $[z^j]\,f(z)^kg(z)^lh(z)$ as $j$, $k$ and $l$ grow to infinity at a comparable rate. Roughly speaking, they show that these coefficients may expose phase transitions (from quadratic-exponential to cubic-exponential) that are determined by subtle asymptotic linear dependences between $j$, $k$ and $l$.

In a different context, Pemantle and Wilson~\cite{PemWil_a,PemWil_b} use residue theory and multidimensional non-coalescing saddle-point methods to analyze the asymptotic behavior of the coefficients of multi-variable generating functions of a meromorphic type. Indeed, since
\[[z^{n_0}]\prod_{j=1}^mf_j^{n_j}=[u^{n_0}v_1^{n_1}\cdots v_m^{n_m}]\,\frac{1}{\prod_{j=1}^m(1-v_j\cdot f_j(u))},\]
under some technical conditions, the asymptotic treatment of the coefficients in~(\ref{ide:problem definition}) is adequate to Pemantle and Wilson's scheme but provided that all the exponents $n_0,\ldots,n_m$ grow to infinity at a comparable rate i.e. $\|\bfn\|=O(n_j)$, for all $0\le j\le m$. In the special case of $m=1$ on the right-hand side above, the partially unpublished dissertation of Lladser~\cite{Lla_a} implies that the condition that $n_0$ and $n_1$ grow to infinity at a comparable rate can be relaxed. This extension of Pemantle and Wilson's scheme for two-variable generating functions follows from an adaptation of the coalescing saddle-point method (see~\cite{Lla_b,Lla_c}).

To the best of our knowledge, the most general discussion of the asymptotic behavior of the coefficients in (\ref{ide:problem definition}) is due to Gardy~\cite{Gar}. However the machinery developed in~\cite{Gar} is
restricted to terms $f_j$ with non-negative Taylor coefficients. Specifically, under the additional assumption that $f_1(0)\cdot f_1'(0)\ne0$ and that $n_j=o(n_1/\sqrt{n_0})$ for $j>1$, Gardy determines the leading order asymptotic term of the coefficients in (\ref{ide:problem definition}) for two different regimes, namely $n_0=\Theta(n_1)$ and $n_0=o(n_1)$. The role of these hypotheses is primarily technical: they prevent the saddle-point of the integral obtained by applying Cauchy's formula in (\ref{ide:problem definition}) to stay away from the saddle-point of integral (also obtained from Cauchy's formula) that represents the coefficients $\prod_{j>1}f_j(\rho)^{n_j}\cdot [z^{n_0}]f_1^{n_1}$. Here $\rho$ is uniquely defined by the relation $\rho f_1'(\rho)/f_1(\rho)=n_0/n_1$.

In this manuscript we introduce a framework that permits a systematic treatment for analyzing the asymptotic behavior of the coefficients in~(\ref{ide:problem definition}) as $\|\bfn\|\to\infty$. The special cases considered in the literature for coefficients of this sort are developed in more generality and without restricting the terms $f_j$ to have non-negative coefficients nor the exponents $n_0,\ldots,n_m$ to be of the same asymptotic order. Within this framework, the asymptotic analysis of the coefficients in~(\ref{ide:problem definition}) is amenable to either the non-coalescing or the coalescing saddle-point method. 

Our main result, Theorem~\ref{thm:parameter integral}, states that it is often possible to relate the asymptotic behavior of the coefficients in (\ref{ide:problem definition}) with the asymptotic behavior of an oscillatory integral of the form
\[\int \exp\left\{-\|\bfn\|\cdot F\left(\theta;\frac{\bfn}{\|\bfn\|}\right)\right\}d\theta,\]
where the integration occurs over a subinterval of $[-\pi,\pi]$ that is centered at the origin.  Theorem~\ref{thm:properties F},
the accompanying result to Theorem~\ref{thm:parameter integral}, states that the function $F$, that we refer to as the {\it associated phase term}, is a continuous function of its arguments however it is also an analytic function of $\theta$. Furthermore, since
$\frac{\partial F}{\partial \theta}\left(0;\bfn/\|\bfn\|\right)=0$ and the $\Re\left\{F\left(\theta;\bfn/\|\bfn\|\right)\right\}$ is minimized at $\theta=0$ the asymptotic behavior of the above integral is amenable to saddle-point methods.

Our methodology was partially developed in~\cite{Lla_a} for a special cases of $m=1$ and $m=2$ in (\ref{ide:problem definition}). The proof of our two main results are omitted but are included in the full version of this manuscript. Furthermore, to keep the present discussion accessible for a wide audience, the application discussed
in Section~\ref{sec:application} relies exclusively on the non-coalescing saddle-point method. The coalescing saddle-point method is relevant only for Section~\ref{sec:planar maps} however we do not implement it in there and rely on known results already available in the literature. We refer the interested reader to~\cite{BFSS_a,BFSS_b} and to Section 5.5 in~\cite{Lla_a} to complement the discussion in Section~\ref{sec:planar maps}.

In what follows it is assumed that for all $1\le j\le m$,
\begin{itemize}
\item[(H1)] $f_j(0)\ne0$, and
\item[(H2)] $f_j$ is a non-constant function of $z$.
\end{itemize}
It should be clear that condition (H2) does not reduce the generality of our exposition. On the contrary, while condition (H1) seems to reduce the generality of our results, this restriction is only apparent (see Remark~\ref{rmk:fjstohjs}).

Our main result develops from the following definition.

\begin{definition}
Given a vector $\bfd\in\RR_+^{m+1}$, we say that a point $z$ with $|z|<r$ is a {\it critical point associated with the direction} $\bfd$ for $(f_1,\ldots,f_m)$ if
\begin{eqnarray}
\label{def:max} \prod_{j=1}^m|f_j(z)|^{d_j}&=&\max_{x:|x|=|z|}\prod_{j=1}^m|f_j(x)|^{d_j},\\
\label{def:critical} d_0&=&\sum\limits_{j=1}^m d_j\cdot \frac{z f_j'(z)}{f_j(z)}.
\end{eqnarray}
If the maximum in (\ref{def:max}) is achieved only at $x=z$, then we say that $z$ is a {\it strictly minimal critical point associated with the direction} $\bfd$.
\end{definition}

Some technical remarks are worth mentioning regarding the above definition. First, observe that if $z$ is a critical point associated with $\bfd$ for $(f_1,\ldots,f_m)$ then it is also a critical point associated with $t\cdot\bfd$, for any $t>0$. Therefore, without loss of generality, we may always assume that $\bfd\in\mathbb{S}_+^m$. In particular, due to condition (H1), it follows that $z=0$ is a strictly minimal critical point associated with any direction $\bfd\in\mathbb{S}_+^m$ such that $d_0=0$. Secondly, observe that due to the maximum modulus principle, (H2) and (\ref{def:max}) imply that
\begin{equation}
\prod_{j=1}^m f_j(z)\ne0,
\end{equation}
whenever $z$ is a critical point associated with some direction $\bfd$ for $(f_1,\ldots,f_m)$. 

Finally, we remark that the occurrence of critical points and more specifically strictly minimal points is somewhat common in the context of combinatorial or probabilistic generating functions. In this case, the functions $f_j$ have nonnegative coefficients implying that condition (\ref{def:max}) is certainly satisfied for each $z\in[0,r)$. In particular, for all nonzero $(d_1,\ldots,d_m)\in\RR_+^m$, each $z\in[0,r)$ turns out to be a critical point associated with the direction
\[\bfd(z):=\frac{(d_0(z),d_1,\ldots,d_m)}{\|(d_0(z),d_1,\ldots,d_m)\|},\]
where
\[d_0(z):=\sum\limits_{j=1}^m d_j\cdot\frac{zf_j'(z)}{f_j(z)}.\]
Furthermore, if there exists $j$ such that $d_j>0$ and $f_j$ is aperiodic, or if there are $j$ and $k$ such that $d_j,d_k>0$ and the period of $f_j$ is relative prime with the period of $f_k$, then each $z\in[0,r)$ turns out to be a strictly minimal critical point associated with the direction $\bfd(z)$ for $(f_1,\ldots,f_m)$.

\section{Main results.}
\label{sec:main}

Our main result is the following one.

\begin{theorem}\label{thm:parameter integral}
{\bf (Parameter-varying integral representation.)} Suppose that conditions (H1) and (H2) are satisfied, that $\D\subset\mathbb{S}_+^m$ is a compact set and that $\Z:\D\to[|z|<r]$ is a continuous function such that for all $\bfd\in\D$, $\Z(\bfd)$ is a strictly minimal critical point associated with the direction $\bfd$ for $(f_1,\ldots,f_m)$. Define
\begin{equation}
\label{thm:def:F}
F(\theta;\bfd):=i\cdot\theta\cdot d_0-\sum_{j=1}^md_j\cdot\ln\left\{\frac{f_j\big(\Z(\bfd)\cdot e^{i\theta}\big)}{f_j\big(\Z(\bfd)\big)}\right\},
\end{equation}
and consider the shortcut notation
\begin{equation}
\label{thm:def:Z}
\Z:=\Z\left(\frac{\bfn}{\|\bfn\|}\right).
\end{equation}
There exists $\epsilon>0$ sufficiently small such that
\begin{itemize}
\item[(a)] {\bf (large exponent)} for all compact sets $\D_1\subset\D\setminus\Z^{-1}\{0\}$ there exists a constant $c>0$ such that
\begin{equation}
\label{thm:large exponent}
[z^{n_0}]\prod_{j=1}^mf_j^{n_j}\!=\!
\frac{\Z^{-n_0}}{2\pi}\prod_{j=1}^m\{f_j(\Z)\}^{n_j}\!\cdot\!\left(\int_{-\epsilon}^\epsilon\exp\left\{-\|\bfn\|\cdot F\left(\theta;\frac{\bfn}{\|\bfn\|}\right)\right\}d\theta+O(e^{-\|\bfn\|\cdot c})\right),
\end{equation}
for all $\bfn\in\|\bfn\|\cdot\D_1$,
\item[(b)] {\bf (small exponent)} for all set $\D_0\subset\big(\Z^{-1}\{0\}+\BB_\epsilon\big)\cap\big(\D\setminus\Z^{-1}\{0\}\big)$,
\begin{equation}
\label{thm:small exponent}
[z^{n_0}]\prod_{j=1}^mf_j^{n_j}=\frac{\Z^{-n_0}}{2\pi}\prod_{j=1}^m\{f_j(\Z)\}^{n_j}\cdot\int_{-\pi}^\pi\exp\left\{-\|\bfn\|\cdot F\left(\theta;\frac{\bfn}{\|\bfn\|}\right)\right\}d\theta,
\end{equation}
for all $\bfn\in\|\bfn\|\cdot\D_0$.
\end{itemize}
\end{theorem}

The terminology of  {\it large} versus {\it small exponent} is motivated by the applications of Theorem~\ref{thm:parameter integral} in the combinatorial or probabilistic setting i.e. when each $f_j$ has nonnegative coefficients. In this
case, under the hypotheses of the theorem, we must have that $\Z(\D)\subset[0,r)$. Suppose that
$\bfn\in\|\bfn\|\cdot\D$ as $\|\bfn\|\to\infty$; in particular, if $\Z$ is as defined in~(\ref{thm:def:Z}) then
\begin{equation}\label{ide:no/n}
\frac{n_0}{\|\bfn\|}=\Z\cdot\sum_{j=1}^m\frac{n_j}{\|\bfn\|}\cdot\frac{f_j'\big(\Z\big)}{f_j\big(\Z\big)}.
\end{equation}
If $\Z$ is confined to a compact set of $(0,r)$, as $\|\bfn\|\to\infty$, then asymptotics for $[z^{n_0}]\prod_{j=1}^mf_j^{n_j}$ fall in the context of part (a) in Theorem~\ref{thm:parameter integral}. However, due to condition (H2), the right-hand side of (\ref{ide:no/n}) vanishes only at $\Z=0$ unless $n_1=\ldots=n_m=0$. In either case, it follows that $n_0$ and $\|\bfn\|$ are of the same order. Therefore, part (a) describes the coefficient of a relatively large power of $z$ in $\prod_{j=1}^mf_j^{n_j}$. On the contrary, if $\Z\to0$ as $\|\bfn\|\to\infty$ then asymptotics for $[z^{n_0}]\prod_{j=1}^mf_j^{n_j}$ fall in the context of part (b). However, in this case, (\ref{ide:no/n}) implies that $n_0=o(\|\bfn\|)$ i.e. part (b) provides information about a relatively small power of $z$ in $\prod_{j=1}^mf_j^{n_j}$.

\begin{remark}\label{rmk:fjstohjs} 
Condition (H1) in Theorem~\ref{thm:parameter integral} is not as restrictive as it may seem. Indeed, if $k_j$ is the degree of vanishing of $f_j$ about $z=0$ then one may rewrite $f_j(z)=z^{k_j}\cdot g_j(z)$, with $g_j(0)\ne0$. As a result,
\begin{equation}
\label{ide:fstohs}
[z^{n_0}]\prod_{j=1}^mf_j^{n_j}=[z^{n_0-\sum_{j=1}^mk_jn_j}]\prod_{j=1}^mg_j^{n_j}.
\end{equation}
The asymptotic analysis of the coefficients on the right-hand side above is now amenable for Theorem~\ref{thm:parameter integral}. Indeed, for all $\bfd\in\mathbb{S}_+^m$, $z$ is a strictly minimal critical point associated with $\bfd$ for $(g_1,\ldots,g_m)$ if and only if $z$ is a strictly minimal critical point associated with $(d_0+\sum_{j=1}^kk_jd_j,d_1,\ldots,d_m)$ for $(f_1,\ldots,f_m)$.

\end{remark}

\begin{remark}\label{rmk:fo}
In some situations one might be more inclined to study the asymptotic behavior of coefficients of the form $[z^{n_0}]\,f_0(z)\cdot\prod_{j=1}^mf_j^{n_j}$ as $\|\bfn\|\to\infty$, where $f_0$ is certain analytic function defined on the disk $[|z|<r]$. In this particular setting, if one defines
\begin{equation}
\label{rmk:amplitude}
A(\theta;\bfd):=f_0(\Z(\bfd)\cdot e^{i\theta}),
\end{equation}
then the asymptotic formulae in (\ref{thm:large exponent}) and (\ref{thm:small exponent}) are respectively replaced by
\begin{eqnarray}
\label{rmk:large exponent}
[z^{n_0}]\,f_0\prod_{j=1}^mf_j^{n_j}\!\!\!\!&=&\!\!\!\!\frac{\Z^{-n_0}}{2\pi}\prod_{j=1}^m\{f_j(\Z)\}^{n_j}\!\!\cdot\!\!\left(\int_{-\epsilon}^\epsilon e^{-\|\bfn\|\cdot F\left(\theta;\frac{\bfn}{\|\bfn\|}\right)}\cdot A\!\!\left(\!\theta;\!\frac{\bfn}{\|\bfn\|}\right)d\theta+O(e^{-\|\bfn\|\cdot c})\right)\!,\\
\label{rmk:small exponent}
[z^{n_0}]\,f_0\prod_{j=1}^mf_j^{n_j}\!\!\!\!&=&\!\!\!\!\frac{\Z^{-n_0}}{2\pi}\prod_{j=1}^m\{f_j(\Z)\}^{n_j}\!\!\cdot\!\!\int_{-\pi}^\pi e^{-\|\bfn\|\cdot F\left(\theta;\frac{\bfn}{\|\bfn\|}\right)}\cdot A\!\!\left(\!\theta;\!\frac{\bfn}{\|\bfn\|}\right)d\theta.
\end{eqnarray}
The above approach could be preferred, for instance, if some of the exponents $n_j$ in (\ref{ide:problem definition}) remain constant as $\|\bfn\|\to\infty$. In this case, $f_0(z)$ would collect all factors associated with a constant exponent. An advantage to this approach stems from the fact that it might be slightly simpler to find strictly minimal points associated with directions in $\mathbb{S}_+^m$ for $(f_1,\ldots,f_m)$ than in $\mathbb{S}_+^{m+1}$ for $(f_0,f_1,\ldots,f_m)$. 
\end{remark}

\begin{remark}
It should be noted that formula (\ref{thm:small exponent}) is not an asymptotic expansion but an equality based on the integral representation of the coefficient $[z^{n_0}]\prod_{j=1}^m f_j^{n_j}$ by means of Cauchy's formula.
\end{remark}

Theorem~\ref{thm:parameter integral} would be of no use without aditional quantitative information about the phase term $F(\theta;\bfd)$. The following accompanying result reveals the principal ingredients required to determine uniform asymptotic expansions (via non-coalescing or coalescing saddle-point methods) for the integral terms appearing in (\ref{thm:large exponent}), (\ref{thm:small exponent}), (\ref{rmk:large exponent}) and (\ref{rmk:small exponent}).

\begin{theorem}
\label{thm:properties F}
{\bf (Associated phase term properties.)} Under conditions (H1) and (H2) it follows that for all open disk $\Theta\subset\CC$ centered at the origin there exists an open disk $\B\subset\Theta$ centered at the origin such that $F(\theta;\bfd)$, as defined in (\ref{thm:def:F}), is a continuous function over the set
\[\Lambda:=\big(\Theta\times\Z^{-1}\B\big)\cup\big(\B\times\D\big).\]
Furthermore,
\begin{itemize}
\item[(a)] for $\bfd\in\D$, $F(\theta;\bfd)$ is an analytic function of $\theta$ such that $F(0;\bfd)=\frac{\partial F}{\partial\theta}(0;\bfd)=0$,
\item[(b)] for $(\theta;\bfd)\in\Lambda$ such that $\theta\in\RR\setminus\{0\}$ and $\bfd\notin\Z^{-1}\{0\}$, the $\Re\{F(\theta;\bfd)\}>0$,
\item[(c)] $F(\theta;\bfd)=\Z(\bfd)\cdot G(\theta;\bfd)$, for a certain continuous function $G(\theta;\bfd)$ defined over $\Lambda$, and
\item[(d)] for $(\theta;\bfd)\in\Lambda$ such that $\bfd\in\Z^{-1}\{0\}$, $G(\theta;\bfd)=\frac{\partial H}{\partial z}(\theta,0;\bfd)-\theta\cdot\frac{\partial^2 H}{\partial z\partial\theta}(0,0;\bfd)$, where
\[H(\theta,z;\bfd):=i\cdot\theta\cdot d_0-\sum_{j=1}^md_j\cdot\ln\left\{\frac{f_j(z\cdot e^{i\theta})}{f_j(z)}\right\}.\]
\end{itemize}
\end{theorem}

\section{Application: Coefficients of a tri-variate generating function}
\label{sec:application}

From a set of cardinality $k$ a collection of $n$ disjoint pairs is named. Let $c(n,k,t)$ be the number of subsets of cardinality $t$ that fail to contain all of the pairs as subsets. By definition $c(n,k,t)=0$ for $2n>k$. On the other hand, a simple inclusion-exclusion argument shows for $2n\le k$ that
\[c(n,k,t)=\sum_{i=0}^n(-1)^i{n\choose i}{k-2i\choose t-2i},\]
where it is understood that ${i\choose j}=0$ unless $0\le j\le i$. Following the lines of~\cite{PemWil_c} in Section 4.10, observe that
\[\sum_{n,i,k,t:2n\le k}{n\choose i}{k-2i\choose t-2i}x^ny^kz^tw^i=\frac{1}{1-y(1+z)}\cdot\frac{1}{1-xy^2(1+2z+(1+w)z^2)}.\]
In particular, by letting $w=-1$, we see that
\begin{equation}
\label{app:def:C}
C(x,y,z):=\sum_{n,k,t:2n\le k}c(n,k,t)x^ny^kz^t=\frac{1}{1-y(1+z)}\cdot\frac{1}{1-xy^2(1+2z)}.
\end{equation}

Up to notational modifications, the coefficients of $C(x,y,z)$ are analyzed in~\cite{PemWil_c} for some special cases where $n$, $k$ and $t$ blow up to infinity at a comparable rate. These coefficients are also analyzed in~\cite{LLPSSW} for a special case where $n$ and $k$ blow up to infinity at a comparable rate, however, $t$ grows sub-linearly with $k$. The discussion that follows extends the asymptotic analysis in~\cite{PemWil_c} and~\cite{LLPSSW}.

The connection between the coefficients of $C(x,y,z)$ and the discussion in the previous two sections is revealed in the following lemma. Its proof can be found in~\cite{LLPSSW} and follows by two consecutive applications of the geometric series.

\begin{lemma}
\label{lem:trivariate to polynomial}
If $2n\le k$, then $[x^ny^kz^t]\,C(x,y,z)=[z^t]\,(1+z)^{k-2n}(1+2z)^n$.
\end{lemma}

In what follows, it is assumed that $2n\le k$. In $\RR^3$ consider the norm $\|\bfd\|:=|d_0|+|d_1|+2|d_2|$, where $\bfd=(d_0,d_1,d_2)$; in particular, $\mathbb{S}_+^2:=\{\bfd\in\RR_+^3:d_0+d_1+2d_2=1\}$. Define $\D$ to be the set of elements $\bfd\in\mathbb{S}_+^2$
for which there exists a strictly minimal critical point associated with $\bfd$ for $(1+z,1+2z)$. An explicit description of $\D$ is revealed by the following result.

\begin{lemma}\label{lem:zd}
$\D=\{\bfd\in\mathbb{S}_+^2:d_1+d_2>d_0\}$ and, the transformation $\Z:\D\to\RR_+$
defined as
\begin{equation}
\label{lem:Z(d)}
\Z(\bfd):=\frac{2d_0}{d_1+2d_2-3d_0+\sqrt{(d_1+2d_2-3d_0)^2+8d_0(d_1+d_2-d_0)}},
\end{equation}
is such that for all $\bfd\in\D$, $\Z(\bfd)$ is a strictly minimal critical point associated with $\bfd$ for $(1+z,1+2z)$. Furthermore,
\begin{equation}
\label{lem:Z-1(0)}
\Z^{-1}\{0\}=\{\bfd\in\RR_+^3:d_0=0\hbox{ and }d_1+2d_2=1\}.
\end{equation}
\end{lemma}

\begin{proof}
Since $f_1(z):=1+z$ and $f_2(z):=1+2z$ are aperiodic polynomials with non-negative coefficients, strictly minimal critical points associated with any direction for $(f_1,f_2)$ must lie in $\RR_+$. In particular, for a fixed $\bfd\in\mathbb{S}_+^2$, $z$ is a strictly minimal critical point associated with $\bfd$ for $(f_1,f_2)$, provided that it is a solution of the equation
\begin{equation}
\label{prf:lem:equation z}
d_0=\frac{d_1z}{1+z}+\frac{2d_2z}{1+2z},\quad z\ge0.
\end{equation}

We claim that the above equation has at most one solution. To show this, and since the above equation has no solution if $\bfd\in\mathbb{S}_+^2$ is such that $d_1=d_2=0$, we may assume without loss of generality that $d_1\ne0$ or $d_2\ne0$. Given any such $\bfd\in\mathbb{S}_+^2$, suppose that $z\ge0$ is a solution of (\ref{prf:lem:equation z}). Then, since $z/(1+z)$ and $2z/(1+2z)$ are strictly increasing functions of $z\ge0$, it follows for $0\le z_1<z<z_2$ that
\[\frac{d_1 z_1}{1+z_1}+\frac{2d_2 z_1}{1+2z_1}<\frac{d_1 z}{1+z}+\frac{2d_2 z}{1+2z}=d_0<\frac{d_1 z_2}{1+z_2}+\frac{2d_2 z_2}{1+2z_2}.\]
This shows our first claim. 

Next, we claim that $\D\subset\{\bfd\in\mathbb{S}_+^2:d_1+d_2>d_0\}$. We show this by contradiction. Indeed, if $d_1+d_2\le d_0$ and (\ref{prf:lem:equation z}) had a solution $z$ then
\[d_1+d_2\le d_0=\frac{d_1 z}{1+z}+\frac{2d_2 z}{1+2z}.\]
In particular,
\[\frac{d_1}{1+z}+\frac{d_2}{1+2z}\le0.\]
Since $\bfd\in\mathbb{S}_+^2$, the above implies that $d_0=1$ and $d_1=d_2=0$. This is inconsistent with 
~(\ref{prf:lem:equation z}) which would imply that $d_0=0$. Therefore, $\D\subset\{\bfd\in\mathbb{S}_+^2:d_1+d_2>d_0\}$ as claimed.

Finally, we show that $\{\bfd\in\mathbb{S}_+^2:d_1+d_2>d_0\}\subset\D$. As a matter of fact, observe that the solutions of~(\ref{prf:lem:equation z}) coincide with the solutions of
\[2(d_1+d_2-d_0)z^2+(d_1+2d_2-3d_0)z-d_0=0,\quad z\ge0.\]
Since for $d_1+d_2>d_0$,
\begin{equation}
\label{prf:z}
z=\frac{\sqrt{(d_1+2d_2-3d_0)^2+8d_0(d_1+d_2-d_0)}-(d_1+2d_2-3d_0)}{4(d_1+d_2-d_0)}
\end{equation}
is a solution to the above equation, the claim follows. This in turn implies that
$\D=\{\bfd\in\mathbb{S}_+^2:d_1+d_2>d_0\}$. Furthermore, since equation (\ref{prf:lem:equation z}) has at most one solution, it follows that $z$, as defined in (\ref{prf:z}), is a a strictly minimal critical point associated with $\bfd$ for $(f_1,f_2)$. The lemma follows by noticing that the conjugate expression of (\ref{prf:z}) coincides with the definition of $\Z(\bfd)$ in (\ref{lem:Z(d)}).
\end{proof}

As a first application of Theorem~\ref{thm:parameter integral} consider a fixed compact set $\K\subset\D\setminus\Z^{-1}\{0\}$. According to Lemma~\ref{lem:trivariate to polynomial}, if $2n\le t$ then $[x^ny^kz^t]\,C(x,y,z)=[z^t]\,(1+z)^{k-2n}(1+2z)^n$. In particular, using Lemma~\ref{lem:zd} and part (a) in Theorem~\ref{thm:parameter integral}, it follows that for all sufficiently small $\epsilon>0$ there exists a constant $c>0$ such that
\begin{equation}
\label{apl:large exponent}
[x^ny^kz^t]\,C(x,y,z)=\Z^{-t}(1+\Z)^{k-2n}(1+2\Z)^n\cdot\Big(I_\epsilon(n,k,t)+O(e^{-\|(t,k-2n,n)\|\cdot c})\Big),
\end{equation}
for all $(n,k,t)$ such that $2n\le k$ and $(t,k-2n,n)/\|(t,k-2n,n)\|\in\K$ as $\|(t,k-2n,n)\|\to\infty$, where
\begin{eqnarray}
\label{apl:cor:def:Z}
\Z&:=&\Z\left(\frac{(t,k-2n,n)}{\|(t,k-2n,n)\|}\right)=\frac{2t}{k-3t+\sqrt{(k-3t)^2+8t(k-n-t)}},\\
\label{apl:integral}
I_\epsilon(n,k,t)&:=&\frac{1}{2\pi}\int_{-\epsilon}^\epsilon\exp\left\{-\|(t,k-2n,n)\|\cdot F\left(\theta;\frac{(t,k-2n,n)}{\|(t,k-2n,n)\|}\right)\right\}d\theta,\\
\label{apl:def:F}
F(\theta;\bfd)&:=&i\cdot d_0\cdot\theta-d_1\cdot\ln\left\{\frac{1+\Z(\bfd) e^{i\theta}}{1+\Z(\bfd)}\right\}-d_2\cdot\ln\left\{\frac{1+2\Z(\bfd) e^{i\theta}}{1+2\Z(\bfd)}\right\},\\
\label{apl:series F}
&=&
\frac{\Z(\bfd)}{2}\left\{\frac{d_1}{(1+\Z(\bfd))^2}+\frac{2d_2}{(1+2\Z(\bfd))^2}\right\}\theta^2+\ldots
\end{eqnarray}
Using the quantitative properties stated for $F(\theta;\bfd)$ in Theorem~\ref{thm:properties F}, one may analyze the asymptotic behavior of $I_\epsilon(n,k,t)$ by means of the non-coalescing saddle-point method. Our findings are summarized in the following result.

\begin{corollary}
\label{cor:large exponent}
Let $\delta\in(0,1)$. If $(n,k,t)$ is such that $2n\le k$, $(1+\delta)t+n\le(1-\delta)k$, and $k=O(t)$ as $t\to\infty$ then
\begin{equation}
\label{apl:cor:large exponent}
[x^ny^kz^t]\,C=\frac{\Z^{-(t+1/2)}(1+\Z)^{k-2n}(1+2\Z)^n}{\sqrt{2\pi}}\cdot\left\{\frac{k-2n}{(1+\Z)^2}+\frac{2n}{(1+2\Z)^2}\right\}^{-1/2}\!\!\!\!\!\!\cdot(1+o(1)).
\end{equation}
\end{corollary}

\begin{proof}
Observe that $(1+\delta)t+n\le(1-\delta)k$ is equivalent to $t+\delta\|(t,k-2n,n)\|\le (k-2n)+n$. As a result, if $t\ne0$ and $\alpha>0$ is such that $k\le\alpha t$ then $(t,k-2n,n)/\|(t,k-2n,n)\|\in\K$, where $\K:=\{\bfd\in\mathbb{S}_+^2:1/(1+\alpha)\le d_0\le (d_1+d_2)-\delta\}$ is a compact set. On the other hand, Lemma~\ref{lem:zd} implies that $\K\subset\D\setminus\Z^{-1}\{0\}$. Since $\|(t,k-2n,n)\|=(k+t)\to\infty$ as $t\to\infty$, formula (\ref{apl:large exponent}) applies for the regime of Corollary~\ref{cor:large exponent} if $\alpha$ is the multiplicative constant in $k=O(t)$.

According to Theorem~\ref{thm:properties F} there exists an open disk $\B\subset\CC$ centered at the origin such that $F(\theta;\bfd)$ is continuous for $(\theta;\bfd)\in\B\times\K$ and analytic for $\theta\in\B$. Without loss of generality assume that $\epsilon$ in (\ref{apl:large exponent}) is such that $[-\epsilon,\epsilon]\subset\B$. Since, for $\bfd\in\K$, $\Z(\bfd)\ne0$,
we see from (\ref{apl:series F}) that 
\begin{equation}\label{apl:partial F at theta=0}
\frac{1}{2!}\frac{\partial^2 F}{\partial\theta^2}(0;\bfd)=\frac{\Z(\bfd)}{2}\left\{\frac{d_1}{(1+\Z(\bfd))^2}+\frac{2d_2}{(1+2\Z(\bfd))^2}\right\}.
\end{equation}

Since the above quantity is strictly positive for $\bfd\in\K$, $F(\theta;\bfd)$ vanishes to constant degree 2 about $\theta=0$ for all $\bfd\in\K$. This motivates us to consider the transformation
\[\tau=\phi(\theta;\bfd):=\theta\cdot\left\{\frac{1}{2!}\frac{\partial^2 F}{\partial\theta^2}(0;\bfd)\right\}^{1/2}\cdot\left(1+\frac{F(\theta;\bfd)-\frac{1}{2!}\frac{\partial^2 F}{\partial\theta^2}(0;\bfd)\cdot\theta^2}{\frac{1}{2!}\frac{\partial^2 F}{\partial\theta^2}(0;\bfd)\cdot\theta^2}\right)^{1/2},\]
where the square-roots are to be interpreted in the principal sense. The above transformation is a well-defined continuous function of $(\theta;\bfd)\in\B\times\K$ and analytic function of $\theta\in\B$ provided that the disk $\B$ is chosen with a sufficiently small radius. Indeed, for a sufficiently small radius, one can warranty that for all $\bfd\in\K$, $\phi(\theta;\bfd)$ is a 1-to-1 analytic function of $\theta\in\B$. Furthermore, since $\phi(0;\bfd)=0$ and, using that $\K$ is compact,
\[\inf_{\bfd\in\K}\left|\frac{\partial\phi}{\partial\theta}(0;\bfd)\right|>0,\]
K{\"o}be 1/4-Theorem implies that there exists a radius $\rho>0$ such that, for all $\bfd\in\K$, $[|\tau|<\rho]\subset\phi(\B;\bfd)$. On the other hand, since for each $\bfd\in\K$, $\phi^{-1}(\tau;\bfd)$ is a  1-to-1 analytic function of $\tau$ such that $\phi^{-1}(0;\bfd)=0$ and
\[\inf_{\bfd\in\K}\left|\frac{\partial\phi^{-1}}{\partial\tau}(0;\bfd)\right|=\frac{1}{\sup\limits_{\bfd\in\K}\left|\frac{\partial\phi}{\partial\theta}(0;\bfd)\right|}>0,\]
K{\"o}be 1/4-Theorem implies that there exists a radius $\delta>0$ such that, for all $\bfd\in\K$, $[|\theta|<\delta]\subset\phi^{-1}([|\tau|<\rho];\bfd)$. Summarizing, we have shown that there exist $\delta,\rho>0$ such that, for all $\bfd\in\K$, $\phi([|\theta|<\delta];\bfd)\subset[|\tau|<\rho]\subset\phi(\B;\bfd)$. Assume without loss of generality that $\epsilon<\delta$. In particular, if $(n,k,t)$ is such that $\bfd=(t,k-2n,n)/\|(t,k-2n,n)\|\in\K$ then, by substituting $\tau=\phi(\theta;\bfd)$ in (\ref{apl:integral}), it follows that
\begin{equation}\label{apl:I_epsilon using tau}
I_\epsilon(n,k,t)=\frac{1}{2\pi}\int_{\phi(-\epsilon;\bfd)}^{\phi(\epsilon;\bfd)}e^{-(k+t)\cdot\tau^2}\,\frac{\partial}{\partial\tau}\big[\phi^{-1}(\tau;\bfd)\big]\,d\tau,
\end{equation}
uniformly for $(t,k-2n,n)/\|(t,k-2n,n)\|\in\K$. Furthermore, since $\epsilon<\delta$, the integrant above is an analytic function of $\tau$ for $|\tau|<\rho$ and hence the line integral above is path independent. On the other hand, observe that if $\tau=\phi(\pm\epsilon;\bfd)$ then $\Re\{\tau^2\}=\Re\{F(\pm\epsilon;\bfd)\}$. Since $\K$ is a compact subset of $\D\setminus\Z^{-1}\{0\}$, Theorem~\ref{thm:properties F} implies that
\[\inf_{\bfd\in\K}\Re\{(\phi(\pm\epsilon;\bfd))^2\}=\inf_{\bfd\in\K}\Re\{F(\pm\epsilon;\bfd)\}>0.\]
Consequently, according to the saddle-point method, the asymptotic behavior of the integral in (\ref{apl:I_epsilon using tau}) is, up to an exponentially decreasing term, equivalent to the asymptotic behavior of the same line integral but with the contour of integration replaced by a small real-neighborhood of $\tau=0$. More specifically, since
\[\left.\frac{\partial}{\partial\tau}\big[\phi^{-1}(\tau;\bfd)\big]\right|_{\tau=0}=\left\{\frac{1}{2!}\frac{\partial^2 F}{\partial\theta^2}(0;\bfd)\right\}^{-1/2},\]
it follows that
\begin{eqnarray*}
I_\epsilon(n,k,t)&=&\left\{\frac{1}{2!}\frac{\partial^2 F}{\partial\theta^2}(0;\bfd)\right\}^{-1/2}\cdot\frac{1}{2\pi}\int_{-\infty}^{\infty}e^{-(k+t)\cdot\tau^2}d\tau\cdot(1+o(1)),\\
&=&\left\{\frac{1}{2!}\frac{\partial^2 F}{\partial\theta^2}(0;\bfd)\right\}^{-1/2}\cdot\frac{1}{2}\sqrt{\frac{1}{\pi\cdot(k+t)}}\cdot(1+o(1)).
\end{eqnarray*}
The corollary follows from (\ref{apl:large exponent}) and (\ref{apl:partial F at theta=0}).
\end{proof}

The second application of Theorem~\ref{thm:parameter integral} is concerned with the asymptotic behavior of $[x^ny^kz^t]\,C$ for $(n,k,t)$ such that $(t,k-2n,n)/\|(t,k-2n,n)\|\in\D$ however $(t,k-2n,n)/\|(t,k-2n,n)\|$ approaches to the set $\Z^{-1}\{0\}$ as $\|(t,k-2n,n)\|\to\infty$. In this case, part (b) in Theorem~\ref{thm:parameter integral} and part (c) in Theorem~\ref{thm:properties F} imply that
\begin{equation}
\label{apl:small exponent}
[x^ny^kz^t]\,C=\Z^{-t}(1+\Z)^{k-2n}(1+2\Z)^n\cdot I(n,k,t),
\end{equation}
where
\begin{eqnarray}
\label{apl:I(n,k,t)}
I(n,k,t)&:=&\frac{1}{2\pi}\int_{-\pi}^\pi\exp\left\{-\|(t,k-2n,n)\|\cdot \Z\cdot G\left(\theta;\frac{(t,k-2n,n)}{\|(t,k-2n,n)\|}\right)\right\}d\theta,\\
\label{apl:series G}
G(\theta;\bfd) &=&
\frac{1}{2}\left\{\frac{d_1}{(1+\Z(\bfd))^2}+\frac{2d_2}{(1+2\Z(\bfd))^2}\right\}\theta^2+\ldots
\end{eqnarray}
and $\Z$ is as defined in (\ref{apl:cor:def:Z}). Observe that if $\Z\to0$ as $\|(t,k-2n,n)\|\to\infty$ then the coefficient $[x^ny^kz^t]\,C$ exposes a phase transition when $\|(t,k-2n,n)\|\cdot \Z$ goes from a bounded quantity to an unbounded one. The properties specified for $G(\theta;\bfd)$ in Theorem~\ref{thm:properties F} lets us obtain the following result.

\begin{corollary}
\label{apl:cor:small exponent}
If $(n,k,t)$ is such that $2n\le k$ and $t=O(1)$ as $k\to\infty$ then
\begin{equation}
\label{apl:cor:t bounded}
[x^ny^kz^t]\,C(x,y,z)=\frac{k^t}{t!}\cdot(1+o(1)).
\end{equation}
If $(n,k,t)$ is such that $2n\le k$ and $t=o(k)$ as $k\to\infty$ and $t\to\infty$ then
\begin{equation}
\label{apl:cor:t to infinity}
[x^ny^kz^t]\,C(x,y,z)=\frac{\Z^{-t}(1+\Z)^{k-2n}(1+2\Z)^n}{\sqrt{2\pi t}}\cdot(1+o(1)).
\end{equation}
\end{corollary}

\begin{proof}
Lemma~\ref{lem:zd} implies that $\Z^{-1}\{0\}=\{\bfd\in\mathbb{S}_+^2:d_0=0\hbox{ and }d_1+2d_2=1\}$. Observe that in the context of formulae (\ref{apl:cor:t bounded}) and (\ref{apl:cor:t to infinity}), $t=o(k)$; in particular, since $2n\le k$, $t<k-n$ for all $k$ sufficiently large. Thus, for $k$ large, Lemma~\ref{lem:zd} implies that $(t,k-2n,n)/\|(t,k-2n,n)\|\in\D$. In addition, since $t/\|(t,k-2n,n)\|=t/(t+k)\to0$ as $k\to\infty$, it follows that $(t,k-2n,n)/\|(t,k-2n,n)\|$ approaches the set $\Z^{-1}\{0\}$ as $k\to\infty$. All this implies that the asymptotic formula in (\ref{apl:small exponent}) applies for all $k$ sufficiently large. On the other hand, part (d) in Theorem~\ref{thm:properties F} implies that $G(\theta;\bfd)=(d_1+2d_2)\cdot(1+i\theta-e^{i\theta})=(1+i\theta-e^{i\theta})$, for all $\bfd\in\Z^{-1}\{0\}$. As a result,
\begin{equation}\label{apl:G uniform limit}
\lim_{k\to\infty}G\left(\theta;\frac{(t,k-2n,n)}{\|(t,k-2n,n)\|}\right)=(1+i\theta-e^{i\theta}),
\end{equation}
and the above limit is uniform for all $\theta\in[-\pi,\pi]$. Furthermore, since $t=o(k)$ then, using (\ref{apl:cor:def:Z}), it follows that $\Z=t/k\cdot(1+o(1))$ as $k\to\infty$. As a result, since $\|(t,k-2n,n)\|=(k+t)=k\cdot(1+o(1))$, the condition that $t=o(k)$ implies
\[\|(t,k-2n,n)\|\cdot \Z\cdot G\left(\theta;\frac{(t,k-2n,n)}{\|(t,k-2n,n)\|}\right)=t\cdot(1+i\theta-e^{i\theta})\cdot(1+o(1)),\]
as $k\to\infty$, uniformly for all $\theta\in[-\pi,\pi]$. 

If $t=O(1)$ as $k\to\infty$ then, in particular, $t=o(k)$. As a result, we may use the above asymptotic formula in (\ref{apl:small exponent}) to obtain that
\begin{eqnarray*}
[x^ny^kz^t]\,C
\!\!\!&=&\!\!\!\Z^{-t}(1+\Z)^{k-2n}(1+2\Z)^n\cdot\frac{1}{2\pi}\int_{-\pi}^\pi e^{-t(1+i\theta-e^{i\theta})\cdot(1+o(1))}d\theta,\\
\!\!\!&=&\!\!\!\left(\frac{k\cdot e}{t}\right)^t\cdot\left\{\frac{1}{2\pi}\int_{-\pi}^\pi e^{-t(1+i\theta-e^{i\theta})}d\theta+o(1)\right\}\cdot(1+o(1)),\\
\!\!\!&=&\!\!\!\left(\frac{k\cdot e}{t}\right)^t\cdot\left\{\frac{1}{2\pi i}\int_{|z|=1}\frac{e^{t(z-1)}}{z^{t+1}}dz+o(1)\right\}\cdot(1+o(1)),\\
\!\!\!&=&\!\!\!\left(\frac{k\cdot e}{t}\right)^t\cdot\left\{\frac{1}{t!}\left(\frac{t}{e}\right)^t+o(1)\right\}\cdot(1+o(1)),\\
\!\!\!&=&\!\!\!\frac{k^t}{t!}\cdot(1+o(1)).
\end{eqnarray*}
This shows (\ref{apl:cor:t bounded}).

Finally, we consider the case where $(n,k,t)$ is such that $2n\le k$ and $t=o(k)$ as $k,t\to\infty$. Since $t=o(k)$ then $\|(t,k-2n,n)\|\cdot\Z=t\cdot(1+t/k)\cdot(1+o(1))\ge t(1+o(1))$; in particular, $\|(t,k-2n,n)\|\cdot\Z\to\infty$. On the other hand, using (\ref{apl:G uniform limit}), it follows that
\[\lim_{k\to\infty}\Re\left\{G\left(\theta;\frac{(t,k-2n,n)}{\|(t,k-2n,n)\|}\right)\right\}=1-\cos\theta,\]
uniformly for all $\theta\in[-\pi,\pi]$. In particular, under the regime we are considering, it follows that for all sufficiently small $\epsilon>0$ there exists a constant $\gamma>0$ such that
\[I(n,k,t)=\frac{1}{2\pi}\int_{-\epsilon}^\epsilon\exp\left\{-(k+t)\cdot \Z\cdot G\left(\theta;\frac{(t,k-2n,n)}{\|(t,k-2n,n)\|}\right)\right\}d\theta+O(e^{-(k+t)\cdot\Z\cdot\gamma}).\]

To analyze the asymptotic behavior of the integral above, observe that from (\ref{apl:series G}) it follows that
\[\frac{1}{2!}\frac{\partial^2G}{\partial\theta^2}(0;\bfd)=\frac{1}{2}\left\{\frac{d_1}{(1+\Z(\bfd))^2}+\frac{2d_2}{(1+2\Z(\bfd))^2}\right\}.\]
In particular, $G(\theta;\bfd)$ vanishes to constant degree 2 about $\theta=0$ for all $\bfd\in\D$ sufficiently close to $\Z^{-1}\{0\}$. By repeating the argument used to show Corollary~\ref{cor:large exponent}, it follows this time that 
\[I(n,k,t)=\left\{\frac{1}{2!}\frac{\partial^2G}{\partial\theta^2}\left(0;\frac{(t,k-2n,n)}{\|(t,k-2n,n)\|}\right)\right\}^{-1/2}\cdot\frac{1}{2}\sqrt{\frac{1}{\pi(k+t)\Z}}\cdot(1+o(1))+O(e^{-(k+t)\cdot\Z\cdot\gamma}).\]
The asymptotic formula in (\ref{apl:cor:t to infinity}) follows from (\ref{apl:small exponent}) and (\ref{apl:G uniform limit}). This completes the proof of the corollary.
\end{proof}

\section{Application: Core size of planar maps}
\label{sec:planar maps}

To end our manuscript we would like to emphasize that the methodology we have provided can also deal with the occurrence of factors with negative Taylor coefficients in (\ref{ide:problem definition}). To fix ideas consider the probability $p_{n,k}$ that a non-separable rooted map with $(n+1)$-edges has a 3-connected non-separable rooted sub-map with $(k+1)$-edges. The asymptotic behavior of coefficients like these is exhaustively studied in~\cite{BFSS_b} as part of a more general apparatus developed by Banderier et al to study the connectivity properties of the root of a random planar map.

Let $M_n$ be the number of non-separable rooted maps with $(n+1)$-edges and $C_k$ the number of 3-connected non-separable rooted maps with $(k+1)$-edges. It follows from Tutte's work that
\[p_{n,k}=\frac{k\cdot C_k}{n\cdot M_n}\cdot[z^{n-1}]\phi(z)^n\psi(z)^{k-1}\psi'(z),\]
where $\phi(z):=(1+z)^3$ and $\psi(z):=z(1-z)$. Since the coefficients $C_k$ and $M_n$ are of a Lagrangean type their asymptotic behavior can be determined using singularity analysis~\cite{FlaOdl}. 
Hence asymptotics for $p_{n,k}$ are reduced to asymptotics for $[z^{n-1}]\phi(z)^n\psi(z)^{k-1}\psi'(z)$ as $n,k\to\infty$. Using the non-coalescing version of the saddle-point method Banderier et al~\cite{BFSS_a} obtain asymptotic formulae for $p_{n,k}$ for the central region $k=n/3$. On the contrary, the coalescing saddle-point method must be used to study the asymptotic behavior of these coefficients in the central region with scaling window of size $n^{2/3}$, namely when $k$ and $n$ satisfy an asymptotic relation of the form
\begin{equation}\label{ide:Airy regime}
k-\frac{n}{3}=O(n^{2/3}).
\end{equation}
For this regime, the leading asymptotic term of the integral obtained by representing $p_{n,k}$ via Cauchy's formula corresponds to the occurrence of two simple-saddles that coalesce at $z=1/2$ as $n\to\infty$. Banderier et al obtain for all finite real numbers $a\le b$ that
\begin{equation}\label{ide:Airy LLT}
\lim_{(n,k)\to\infty}\sup\limits_{a\le\frac{k-n/3}{n^{2/3}}\le b}\left|n^{2/3}\cdot p_{n,k}-\frac{16}{81}\cdot\frac{3^{4/3}}{4}
\cdot\mathcal{A}
\left(\frac{3^{4/3}}{4}\cdot\frac{k-n/3}{n^{2/3}}\right)\right|=0\,,
\end{equation}
where $\mathcal{A}(x):=2\,e^{-2x^3/3}\,\big(x\cdot\hbox{Ai}(x^2)-\hbox{Ai}'(x^2)\big)$ 
is the density of the map-Airy distribution and $\hbox{Ai}(z)$ is the Airy function (see~\cite{BFSS_a} for more details).

The end of our manuscript is devoted to show how the above analysis would fit in our framework. For this observe that
\begin{equation}\label{ide:linear combination}
[z^{n-1}]\phi(z)^n\psi(z)^{k-1}\psi'(z)=[z^{n_0}](1+z)^{n_1}(1-z)^{n_2}-2[z^{n_0-1}](1+z)^{n_1}(1-z)^{n_2},
\end{equation}
where $n_0:=(n-k)$, $f_1(z):=(1+z)$, $n_1:=3n$, $f_2(z):=(1-z)$ and $n_2:=(k-1)$. In what follows we use $\|\cdot\|_\infty$ as the reference norm; in particular, $\mathbb{S}^2_+=\{\bfd\in\RR_+^3:0\le\min\{d_0,d_1,d_2\}\le\max\{d_0,d_1,d_2\}=1\}$. 
Notice that for the regime in (\ref{ide:Airy regime}), $\|(n_0,n_1,n_2)\|_\infty=3n$ and $(n_0,n_1,n_2)\sim n\cdot (2/3,3,1/3)$. In particular, under the restriction imposed by (\ref{ide:Airy regime}) it applies that
\begin{equation}\label{ide:main boundary point}
\lim_{n\to\infty}\frac{(n_0,n_1,n_2)}{\|(n_0,n_1,n_2\|_\infty}=\left(\frac{2}{9},1,\frac{1}{9}\right).
\end{equation}
Furthermore, observe that the above limit belongs $(\D\cap\partial\D)$, where $\D$ is the set defined as
\[\D:=\{\bfd\in\mathbb{S}^2_+:d_1=1\hbox{ and }(d_1-d_2)^2\ge 4d_0(d_1+d_2-d_0)\}\setminus\{(0,1,1),(1,1,0)\}.\]
The link between the Airy phenomena as described by Banderier et al and our framework is revealed by the following result.

\begin{lemma}
\label{lem:airy regime}
The transformation $\Z:\D\to\RR_+$ defined as
\[\Z(\bfd):=\frac{2d_0}{(d_1-d_2)+\sqrt{(d_1-d_2)^2-4d_0(d_1+d_2-d_0)}}\]
is such that for all $\bfd\in\D$, $\Z(\bfd)$ is a strictly minimal critical point associated with $\bfd$ for $(1+z,1-z)$. Furthermore, 
\[\frac{\partial^2 F}{\partial\theta^2}(0;\bfd)\ge0\]
for all $\bfd\in\D$ with equality only for $\bfd$ such that $d_0=0$ or $(d_1-d_2)^2= 4d_0(d_1+d_2-d_0)$.
\end{lemma}
\begin{proof}
Set $f_1(z):=1+z$ and $f_2(z):=1-z$. Consider $\bfd=(d_0,d_1,d_2)\in\D$; in particular, $d_1=1$ and $d_2<1$. Furthermore, $2d_0\le(1-\sqrt{d_2})^2$ or $2d_0\ge(1+\sqrt{d_2})^2$ and if either of these inequalities is indeed an equality then $\bfd\in\partial\D$. In what follows it is assumed that $2d_0\le(1-\sqrt{d_2})^2$. The proof of the lemma for the other case is similar and left to the reader.

In order for $z$ to be a critical point associated with $\bfd$ for $(f_1,f_2)$ it is necessary that $d_0=d_1z/(1+z)-d_2z/(1-z)$. A straightforward calculation shows that $z=\Z(\bfd)$ is a non-negative solution of this equation. On the other hand, by reducing to polar coordinates, one determines for $r\ge0$ that $|f_1(x)|^{d_1}|f_2(x)|^{d_2}$ has $x=r$ as its unique maxima on the circle $|x|=r$ provided that
\begin{equation}\label{ide:key inequality}
\frac{1-d_2}{1+d_2}\ge\frac{2r}{1+r^2}.
\end{equation}
Solving for the above inequality in terms of the variable $r$, one sees that the condition $r\le(1-\sqrt{d_2})/(1+\sqrt{d_2})$ is sufficient for (\ref{ide:key inequality}) and that the above inequality is strict unless $r=(1-\sqrt{d_2})/(1+\sqrt{d_2})$. In particular, since
\begin{eqnarray*}
\Z(\bfd)&\le&\frac{(1-\sqrt{d_2})^2}{1-d_2+\sqrt{(1-d_2)^2-4d_0(1+d_2-d_0)}},\\
&=&\frac{1-\sqrt{d_2}}{1+\sqrt{d_2}}\cdot\frac{1-d_2}{1-d_2+\sqrt{(1-d_2)^2-4d_0(1+d_2-d_0)}},\\
&\le&\frac{1-\sqrt{d_2}}{1+\sqrt{d_2}},
\end{eqnarray*}
it follows that $\Z(\bfd)$ is a strictly minimal simple point associated with $\bfd$ for $(f_1,f_2)$. The above shows in particular that $\Z(\bfd)<1$. In addition, since
\[\frac{\partial^2 F}{\partial\theta^2}(0;\bfd)=\frac{(1-d_2)\Z(\bfd)}{2(1-\Z(\bfd)^2)^2}\left(\Z(\bfd)-\frac{1-\sqrt{d_2}}{1+\sqrt{d_2}}\right)\left(\Z(\bfd)-\frac{1+\sqrt{d_2}}{1-\sqrt{d_2}}\right),\]
it also follows that $\frac{\partial^2 F}{\partial\theta^2}(0;\bfd)\ge0$. Furthermore, $\frac{\partial^2 F}{\partial\theta^2}(0;\bfd)=0$ if and only if $\Z(\bfd)=0$ or $\Z(\bfd)=(1-\sqrt{d_2})/(1+\sqrt{d_2})$. The first condition is equivalent to having $d_0=0$. The second condition is equivalent to having $2d_0=(1-\sqrt{d_2})^2$ which in turn implies that $\bfd\in\partial\D$. This completes the proof.
\end{proof}

With the aid of the above lemma and of theorems~\ref{thm:parameter integral} and~\ref{thm:properties F}, we see from (\ref{ide:linear combination}) that the asymptotic behavior of $[z^{n-1}]\phi(z)^n\psi(z)^{k-1}\psi'(z)$ can be reduced to the one of an oscillatory integral provided that $n\to\infty$ and that for all $n$ sufficiently large,
\[\left(\frac{n-k}{3n},1,\frac{k-1}{3n}\right)\in\D.\]
Indeed, as it follows from the proof of the Lemma~\ref{lem:airy regime}, the above condition is equivalent to request that $1\le k\le n$ and that
\[\frac{2(n-k)}{3n}\le\left(1-\sqrt{\frac{k-1}{3n}}\right)^2\hbox{ or }\,\frac{2(n-k)}{3n}\ge\left(1+\sqrt{\frac{k-1}{3n}}\right)^2.\]

The regime in (\ref{ide:Airy regime}) is a particular instance of the left-hand side above. In terms of our framework, the occurrence of the Airy function in (\ref{ide:Airy LLT}) is due to (\ref{ide:main boundary point}) and the fact that
\[F\left(\theta;\frac{1}{9},1,\frac{2}{9}\right)=\frac{8i}{81}\theta^3+\frac{10}{81}\theta^4+\ldots\]
The local limit in~(\ref{ide:Airy LLT}) results from part (a) in Theorem~\ref{thm:parameter integral} and from a second-order asymptotic expansion (obtained using the coalescing saddle-point method) of the integrals associated with the two terms on the right-hand side of~(\ref{ide:linear combination}). The need for second-order expansions is because the factor $\psi'(z)=1-2z$ on the right-hand side of (\ref{ide:linear combination}) cancels out at $z=\Z(1/9,1,2/9)=1/2$.
 
\bibliographystyle{alpha}
\bibliography{biblio}

\newcommand{\etalchar}[1]{$^{#1}$}
\begin{thebibliography}{LP{\v S}{\etalchar{+}}06}

\bibitem[BFSS00]{BFSS_a}
C.~Banderier, P.~Flajolet, G.~Schaeffer, and M.~Soria.
\newblock Planar maps and {A}iry phenomena.
\newblock In {\em Automata, Languages and Programming}, pages 388--402, 2000.

\bibitem[BFSS01]{BFSS_b}
C.~Banderier, P.~Flajolet, G.~Schaeffer, and M.~Soria.
\newblock Random maps, coalescing saddles, singularity analysis, and airy
  phenomena.
\newblock {\em Random Structures and Algorithms 19(3-4), 194--246}, 2001.

\bibitem[BH86]{BleHan}
N.~Bleistein and R.~Handelsman.
\newblock {\em Asymptotic expansion of integrals}.
\newblock Dover Publications, 1986.

\bibitem[FO90]{FlaOdl}
P.~Flajolet and A.~M. Odlyzko.
\newblock Singularity analysis of generating functions.
\newblock {\em SIAM Journal on Discrete Mathematics}, 3(2):216--240, 1990.

\bibitem[FS06]{FlaSed}
P.~Flajolet and R.~Sedgewick.
\newblock {\em Analytic Combinatorics}, 2006.
\newblock Electronic version available on-line at
  \href{http://algo.inria.fr/flajolet/Publications/book060418.pdf}{\texttt{htt%
p://algo.inria.fr/flajolet/Publications/book060418.pdf}}.

\bibitem[Gar95]{Gar}
D.~Gardy.
\newblock Some results on the asymptotic behavior of coefficients of large
  powers of functions.
\newblock {\em Discrete Mathematics, vol 139, 189--217}, 1995.

\bibitem[GJ04]{GouJac}
I.~P. Goulden and D.~M. Jackson.
\newblock {\em Combinatorial enumeration}.
\newblock Dover Publications, 2004.

\bibitem[JK77]{JohKot}
N.L. Johnson and S.~Kotz.
\newblock {\em Urn models and their applications}.
\newblock Wiley, 1977.

\bibitem[Lla03]{Lla_a}
M.~Lladser.
\newblock Asymptotic enumeration via singularity analysis.
\newblock {\em Doctoral dissertation, Ohio State University}, 2003.

\bibitem[Lla06a]{Lla_b}
M.~Lladser.
\newblock Uniform formulae for coefficients of meromorphic functions in two
  variables. {P}art {I}.
\newblock {\em To appear in SIAM J. Disc. Mathematics}, 2006.

\bibitem[Lla06b]{Lla_c}
M.~Lladser.
\newblock Uniform formulae for coefficients of meromorphic functions in two
  variables, part {II}: the {A}iry phenomena.
\newblock {\em In preparation}, 2006.

\bibitem[LP{\v S}{\etalchar{+}}06]{LLPSSW}
M.~Lladser, P.~Poto{\v c}nik, J.~{\v S}ir\'a{\v n}, J.~{\v S}iagiov\'a, and
  M.~Wilson.
\newblock The diameter of random {C}ayley digraphs of given degree.
\newblock {\em Preprint}, 2006.

\bibitem[PW02]{PemWil_a}
R.~Pemantle and M.~Wilson.
\newblock Asymptotics of multivariate sequences, part {I}: smooth points of the
  singular variety.
\newblock {\em J. Comb. Theory, Series A, vol. 97, 129--161}, 2002.

\bibitem[PW04]{PemWil_b}
R.~Pemantle and M~Wilson.
\newblock Asymptotics of multivariate sequences, part {II}: multiple points of
  the singular variety.
\newblock {\em Combinatorics, Probability and Computing 13, 735-761}, 2004.

\bibitem[PW05]{PemWil_c}
R.~Pemantle and M.~Wilson.
\newblock Twenty combinatorial examples of asymptotics derived from
  multivariate generating functions.
\newblock {\em Preprint}, 2005.

\bibitem[Wil90]{Wil}
Herbert~S. Wilf.
\newblock {\em Generatingfunctionology}.
\newblock Academic Press, 1990.

\end{thebibliography}

\end{document}